\documentstyle[11pt]{article}
\begin{document}
\setlength{\baselineskip}{18pt}
\title{ Some global results on \\ Quasipolynomial discrete systems }
\author{\\ Benito Hern\'{a}ndez--Bermejo$^{\:a,1}$ \and \\ L\'{e}on Brenig$^{\:b,2}$}
\date{}
\maketitle
\begin{flushleft}
{\em $^{a}$ Escuela Superior de Ciencias Experimentales y Tecnolog\'{\i}a. Edificio Departamental II. } \\
{\em \hspace{0.2 cm} Universidad Rey Juan Carlos. Calle Tulip\'{a}n S/N.
28933--M\'{o}stoles--Madrid. Spain.} \\
\mbox{} \\
{\em $^{b}$ Service de Physique Th\'{e}orique et Math\'{e}matique. Universit\'{e} Libre de Bruxelles.} \\
{\em \hspace{0.2 cm} Campus Plaine -- CP 231. Boulevard du Triomphe. B-1050 Bruxelles. Belgium.}
\end{flushleft}

\mbox{}

\noindent \rule{15cm}{0.01in}
\begin{center}
{\bf Abstract}
\end{center}

The quasipolynomial (QP) generalization of Lotka-Volterra discrete-time systems is considered. Use of the QP formalism is made for the investigation of various global dynamical properties of QP discrete-time systems including permanence, attractivity, dissipativity and chaos. The results obtained generalize previously known criteria for discrete Lotka-Volterra models.

\noindent \rule{15cm}{0.01in}

\mbox{}


\vfill

\noindent $^{1}$ Corresponding author. Telephone: (+34) 91 488 73
91. Fax: (+34) 91 488 73 38. \newline \mbox{} \hspace{0.35mm}
E-mail: {\tt benito.hernandez@urjc.es }

\noindent $^{2}$ E-mail: {\tt lbrenig@ulb.ac.be }

\pagebreak
\begin{flushleft}
{\bf 1. Introduction}
\end{flushleft}

The use of Lotka-Volterra (LV in what follows) discrete-time systems is a well-known subject of
applied Mathematics. They were first introduced in a biomathematical context by Moran
\cite{mor1}, and later popularized by May and collaborators \cite{may1,mo1}. Since then LV systems have proved to be a rich source of analysis for the investigation of dynamical properties and modelling in different domains, not only population dynamics
\cite{fg1,bf1,hhj1,lw1,gs1} but also Physics \cite{raj1,ur1}, Chemistry \cite{gb1} and Economy \cite{doh1}.

Recently the family of discrete-time systems termed quasipolynomial (QP from now on) has deserved some attention in the literature \cite{bl1,bl2,bl3}. In this context it is worth noting that the interest of QP discrete systems arises from several different features.
In first place, they constitute a wide generalization of LV models. However, LV systems are not just a particular case of QP ones but play a central, in fact canonical role in the QP framework, as it will be appreciated in what is to follow. Additionally, QP systems arise as the result of a discretization of the QP ordinary differential equations,
\begin{equation}
\label{qpodes}
   \frac{\mbox{d}x_i}{\mbox{d}t} = x_{i} \left( \lambda _{i} + \sum_{j=1}^{m}A_{ij}\prod_{k=
      1}^{n}x_{k}^{B_{jk}} \right) , \;\:\;\: i = 1 \ldots n
\end{equation}
where $n$ and $m$ are positive integers, and $A$, $B$ and $\lambda$ are real matrices. Apart from their close connection with QP discrete systems, QP differential equations (\ref{qpodes}) have been applied in different domains including integrability properties, normal form analysis, stability, Hamiltonian dynamics, neural networks, biochemical modelling, population dynamics, etc. For instance, the reader is referred to \cite{bl1} for a sample of references on the subject. Finally, QP discrete equations display a set of properties (termed QP formalism) that provide a methodology for the establisment and transfer of results in a matrix algebraic framework \cite{bl1,bl2,bl3}. These features constitute the basis of the developments to be presented in this paper and some of them will be illustrated in the next sections.

The purpose of this work is to make use of the QP techniques in order to demonstrate several results regarding the global dynamics of QP discrete-time systems. The structure of the article is as follows. In Section 2 the foundations of the QP formalism for discrete systems are briefly outlined in order to make the work self-contained. Section 3 is devoted to a first illustration of the use of the QP formalism, leading to the establishment of new criteria about  permanence and global attractivity for QP systems of dimension 2. In Section 4 the extension of this kind of results to the $n$-dimensional case leads to the establishment of properties regarding global attractivity and dissipativity. Finally Section 5 explores further consequences of the previous results on permanence and global attractivity that lead to the demonstration of some chaos criteria.

%\mbox{}

%\mbox{}
\pagebreak
\begin{flushleft}
{\bf 2. Overview and notation of the QP formalism for discrete systems}
\end{flushleft}

The aim of this section is to present an overview of the discrete-time version of the QP formalism. The reader is referred to \cite{bl1} for the full details and for additional aspects of the theory not considered here and to \cite{bl2,bl3} for some applications of the formalism. In this work we shall deal with QP systems of the form
\begin{equation}
\label{qpm}
    x_i(t+1)=x_i(t) \exp \left( \lambda _i + \sum _{j=1}^n A_{ij} \prod _{k=1}^n
    [x_k(t)]^{B_{jk}} \right)
    \:\: , \:\:\:\:\: i =1, \ldots , n
\end{equation}
where {\em (i)} index $t$ is an integer denoting the discrete time; {\em (ii)} variables $x_i(t)$ are assumed to be positive for $i=1, \ldots ,n$ and for every $t$; and {\em (iii)} $A=(A_{ij})$, $B=(B_{ij})$ and $\lambda = ( \lambda _i)$ are real matrices of dimensions $n \times n$, $n \times n$ and $n \times 1$, respectively. The terms
\[
    \prod _{k=1}^n [x_k(t)]^{B_{jk}} \: , \:\:\: j=1, \ldots ,n
\]
appearing in the exponential of equation (\ref{qpm}) are known as quasimonomials. More general QP systems with a number of quasimonomials that may be different from the dimension $n$ are considered in \cite{bl1}. However, some of the results to be demonstrated in the forthcoming sections are not valid when the number of quasimonomials differs from $n$. Consequently, in this article we shall focus on the case (\ref{qpm}) in which they are equal. Notice that LV models
\begin{equation}
\label{lvm}
    x_i(t+1)=x_i(t) \exp \left( \lambda _i + \sum _{j=1}^n A_{ij} x_j(t) \right)
    \:\: , \:\:\:\:\: i =1, \ldots , n
\end{equation}
are a particular case of QP discrete systems, namely the one corresponding to $B$ the $n \times n$ identity matrix.

An important basic property is that the positive orthant is an invariant set for every QP system. This is natural in many domains (such as population dynamics or chemical kinetics) in which variables are positive by definition. In the QP context, this feature is always present. Unless otherwise stated, it is always assumed that QP systems are defined in int$\{ I \!\!
R^n_+\}$.

A key set of transformations relating QP systems are the quasimonomial transformations (QMTs) defined as:
\begin{equation}
\label{qmt}
    x_i(t)= \prod_{j=1}^n [y_j(t)]^{C_{ij}} \:\: , \:\:\: i = 1, \ldots ,n \:\: ; \:\:\:
    \mid C \mid \neq 0
\end{equation}
The reader is referred to \cite{bl1} for other types of QP formalism transformations that will not be employed in this work. The form-invariance of QP systems after a QMT (\ref{qmt}) is one of the cornerstones of the formalism. Actually, if we consider a $n$-dimensional QP system of matrices $A$, $B$, $\lambda$ and perform a QMT of matrix $C$, the result is another
$n$-dimensional QP system of matrices $A'$, $B'$, $\lambda '$ where:
\begin{equation}
\label{trmat}
    A' = C^{-1} \cdot A \: , \:\:\:
    B' = B \cdot C \: , \:\:\:
    \lambda ' = C^{-1} \cdot \lambda
\end{equation}
Moreover, every QMT relating two QP systems is a topological conjugacy. Consequently, we not only have a formal invariance between QP systems related by a QMT, but actually a complete dynamical equivalence (in the topological sense). These properties imply that the set of all QP systems related by means of QMTs actually constitute an equivalence class. One important label of such classes is given by the matrix products $\Gamma = B \cdot A$ and $\Lambda = B \cdot
\lambda$, which are invariant for every equivalence class.

Notice that if on a QP system (\ref{qpm}) we perform a QMT of matrix $C=B^{-1}$ (assuming that $B$ is invertible) then the result is a LV model (\ref{lvm}) with matrices $A_{LV}=B \cdot A$ and $\lambda _{LV} = B \cdot \lambda$. Therefore every QP system is equivalent \cite{bl1} to a LV system, where this equivalence is not only formal but complete in a topological sense, as indicated before. Moreover, the LV matrices are precisely the matrix invariants $\Gamma = B \cdot A$ and $\Lambda = B \cdot \lambda$ of the equivalence class. Consequently LV equations are in fact canonical representatives of the QP equivalence classes. This property allows the use of the formalism for the establishment of new results for QP systems, as it will be seen in what follows.

\mbox{}

\mbox{}

\begin{flushleft}
{\bf 3. Permanence and global attractivity in dimension 2}
\end{flushleft}

As a first illustration of the QP operational framework for the derivation of results we can consider the case $n=2$. It is worth recalling that a QP system is said to be permanent if there is a compact set ${\cal S} \subset \mbox{int} \{ I \!\! R^n_+ \}$ such that for every initial condition $x(0) \in \mbox{int} \{ I \!\! R^n_+ \}$ the orbits enter and remain within ${\cal S}$. We recall also the definition of pattern of a $n \times m$ real matrix $P$, which is an array Pattern($P$) of the same dimensions whose entry at the position $(i,j)$ is one of  the symbols $\{ +, -, 0 \}$ depending on the positive, negative or zero value of element
$P_{ij}$, respectively. We start with a first result for the cooperative case:

\mbox{}

\noindent {\bf Theorem 1.} {\em Consider a QP system (\ref{qpm}) with $n=2$ and invertible matrix $B$. If the class invariant $\Gamma$ has the cooperative pattern}
\begin{equation}
\label{coopat}
    \mbox{Pattern} ( \Gamma ) = \left( \begin{array}{cc} - & + \\ + & - \end{array} \right)
\end{equation}
{\em then the system is not permanent.}

\mbox{}

\noindent {\bf Proof.} Since $B$ is invertible, the application to the QP system of a QMT of matrix $C=B^{-1}$ reduces it to a LV system according to the transformation rules
(\ref{trmat}). Now notice that matrix $A_{LV}$ of the resulting LV system is the class invariant $\Gamma$. But for such kind of LV systems it is known that they are not permanent (as far as every two-dimensional discrete-time LV system such that its matrix $A_{LV}$ has pattern (\ref{coopat}) is not permanent \cite{lw1}). Then the topological conjugacy property of QMTs \cite{bl1} mentioned in Section 2 maps the property of not being permanent to the initial QP system and the proof is complete. $\:\:\: \Box$

\mbox{}

We now consider a complementary case, namely the competitive one:

\mbox{}

\noindent {\bf Theorem 2.} {\em Consider a QP system (\ref{qpm}) with $n=2$ and invertible matrix $B$. If the following three conditions hold:

\noindent (a) The class invariant $\Gamma$ has the competitive pattern:
\begin{equation}
    \label{patcom}
    \mbox{\rm Pattern} ( \Gamma ) = \left( \begin{array}{cc} - & - \\ - & - \end{array}     \right)
\end{equation}

\noindent (b) The class invariant $\Lambda$ has positive pattern:
\begin{equation}
    \label{placom}
    \mbox{\rm Pattern} ( \Lambda ) = \left( \begin{array}{c} + \\ + \end{array} \right)
\end{equation}

\noindent (c) The matrix product $\Gamma ^T \cdot \Lambda ^{\bot}$ has pattern
\[
    \mbox{\rm Pattern} ( \Gamma ^T \cdot \Lambda ^{\bot} ) =
    \left( \begin{array}{c} - \\ + \end{array} \right)
\]
where superscript $^T$ denotes the transpose matrix and superscript $^{\bot}$ is defined as
\[
    \Lambda ^{\bot} = \left( \begin{array}{c} \Lambda _1 \\ \Lambda _2 \end{array} \right)
     ^{\bot} = \left( \begin{array}{r} \Lambda _2 \\ - \Lambda _1 \end{array} \right)
\]
Then the system is permanent.}

\mbox{}

\noindent {\bf Proof.} Again the fact that $B$ is invertible allows the application to the QP system of a QMT of matrix $C=B^{-1}$, and according to identities (\ref{trmat}) the outcome will be a LV system. Now the coefficient matrices $A_{LV}$ and $\lambda_{LV}$ of the resulting LV system are the class invariants $\Gamma$ and $\Lambda$, respectively. Note that
two-dimensional LV discrete-time systems are permanent if the following three conditions are satisfied \cite{lw1}: (i) having the competitive pattern (\ref{patcom}) for its matrix $A_{LV}$; (ii) having pattern (\ref{placom}) for its matrix $\lambda_{LV}$; and (iii) verifying the conditions
\begin{equation}
    \label{clvpt}
    \lambda_2 A_{11}-\lambda_1 A_{21}<0 \;\;\: , \;\;\:\: \lambda_1 A_{22}-\lambda_2 A_{12}<0
\end{equation}
Consider now the initial QP system complying to the conditions of Theorem 2. After the QMT of matrix $C = B^{-1}$ the result is a LV system verifying conditions (i) to (iii) above because inequalities (\ref{clvpt}) are equivalent to condition (c) of the theorem. Consequently such LV system is permanent. The topological conjugacy property of QMTs \cite{bl1} maps the permanence property to the initial QP system and the result is demonstrated. $\:\:\: \Box$

\mbox{}

Notice that conditions (a) to (c) of Theorem 2 (and equivalently conditions (i) to (iii) of the Proof for LV systems) imply that there exists a unique fixed point in $\mbox{int} \{ I \!\! R^2_+ \}$. Following a similar strategy to the results above, it is possible to demonstrate also an additional result regarding global attractivity for the same kind of QP systems:

\mbox{}

\noindent {\bf Theorem 3.} {\em Consider a QP system (\ref{qpm}) with $n=2$ and invertible matrix $B$ verifying the same hypotheses (a) to (c) of Theorem 2. If in addition $\Lambda_i \leq 1$ for $i=1,2$, then the unique positive fixed point of the QP system is globally asymptotically stable in $\mbox{int} \{ I \!\! R^2_+ \}$.}

\mbox{}

\noindent {\bf Proof.} The first step is the reduction of the QP system given to its equivalent LV canonical representative by means of a QMT of matrix $C=B^{-1}$. Such two-dimensional LV equations have competitive pattern (\ref{patcom}) for its matrix $A_{LV}$ and pattern
(\ref{placom}) for its matrix $\lambda_{LV}$, and verify conditions (\ref{clvpt}) as well (which are equivalent to condition (c) of Theorem 2 in the LV case). For LV systems verifying these conditions it is well-known that if in addition $\lambda_i \leq 1$ for $i=1,2$ then the unique positive fixed point is a global attractor in $\mbox{int} \{ I \!\! R^2_+ \}$ (see
\cite{wwzl}). This is therefore the case for the LV system obtained after the QMT, as implied by the hypotheses of the theorem. Then, such result is generalized to every QP system verifying the conditions of Theorem 3 due to the QMT topological conjugacy property \cite{bl1}. $\:\:\: \Box$

\mbox{}

It is worth noting that it is possible to imagine situations in which Theorem 2 is applicable but Theorem 3 is not. In such cases it is clear that the permanence is maintained, but the global attractivity is not necessarily preserved, and actually the system dynamics may become chaotic. We shall consider this kind of possibility later in Section 5.

The last result of this section deals with a third scenario, namely the predator-prey case, which complements the ones previously considered:

\mbox{}

\noindent {\bf Theorem 4.} {\em Consider a QP system (\ref{qpm}) with $n=2$, invertible matrix $B$ and such that the following conditions hold:

\noindent (a) The class invariant $\Gamma$ has the predator-prey pattern:
\begin{equation}
    \label{patpp}
    \mbox{\rm Pattern} ( \Gamma ) = \left( \begin{array}{cc} - & - \\ + & - \end{array}     \right)
\end{equation}

\noindent (b) The class invariant $\Lambda$ has pattern:
\begin{equation}
    \label{plapp}
    \mbox{\rm Pattern} ( \Lambda ) = \left( \begin{array}{c} + \\ - \end{array} \right)
\end{equation}

\noindent (c) $\Lambda_1 \leq 1$.

\noindent (d) $\Lambda_1 \Gamma_{21}/ \Gamma_{11} < \Lambda _2 \leq \Gamma_{21}/ \Gamma_{11}
+1$.

\noindent (e) $\Gamma_{11} \Gamma_{22} + \Gamma_{12} \Gamma_{21} >0$.

\noindent Then the QP system has a unique positive fixed point, and such point is globally asymptotically stable in $\mbox{int} \{ I \!\! R^2_+ \}$.}

\mbox{}

\noindent {\bf Proof.} Since $B$ is invertible a QMT of matrix $C = B^{-1}$ can be performed in order to reduce the system to LV form with matrices $A_{LV}= \Gamma$ and $\lambda_{LV}= \Lambda$, where $\Gamma$ and $\Lambda$ are the QP class invariants of the initial QP system. Consequently the LV matrices $A_{LV}$ and $\lambda_{LV}$ have predator-prey patterns (\ref{patpp}) and (\ref{plapp}), respectively. For such kind of LV systems it can be demonstrated after some simple algebra that there exists a unique fixed point. It is also known that the fixed point is positive and globally asymptotically stable in $\mbox{int} \{ I \!\! R^2_+ \}$ for such LV systems provided the following conditions are verified for the LV matrices $A_{LV}$ and $\lambda_{LV}$ \cite{wwzl}: (i) $\lambda_1 \leq 1$; (ii) $\lambda_1 A_{21}/A_{11} < \lambda_2 \leq A_{21}/A_{11}+1$; and (iii) $A_{11}A_{22}+A_{12}A_{21}>0$. Conditions (i) to (iii) are also verified by the LV system obtained after the QMT of matrix $C = B^{-1}$ as a consequence of hypotheses (c) to (e) of the theorem. The proof is concluded by noting that the result is then generalized for every QP system in the class due to the QMT topological conjugacy property \cite{bl1}. $\:\:\: \Box$

\mbox{}

We conclude the section with an illustration of the previous results:

\mbox{}

\noindent {\bf Example 1.} In \cite{lw1} the following permanent competitive LV system is considered:
\[
    A_{LV1} = \left( \begin{array}{cc} -1 & -1/2 \\ -1/2 & -1 \end{array} \right)
    \:\:\: , \:\:\:\:
    \lambda _{LV1}= \left( \begin{array}{c} 1/2 \\ 1/2 \end{array} \right)
\]
This system verifies the conditions of Theorems 2 and 3 and therefore it is permanent and has a
unique positive fixed point $(1/3,1/3)$ which is in addition globally asymptotically stable in
$\mbox{int} \{ I \!\! R^2_+ \}$. Actually the following generalization (still of LV form) is also permanent:
\begin{equation}
\label{lv2inv}
    A _{LV2}= \left( \begin{array}{cc} -1 & -1/2 \\ -1/2 & -1 \end{array} \right)
    \:\:\: , \:\:\:\:
    \lambda _{LV2}= \left( \begin{array}{c} \rho \\ \rho \end{array} \right)
    \:\: , \:\:\: \rho >0
\end{equation}
The QP formalism allows considering more general systems, i.e. without the restriction of being of LV form. This is equivalent to allowing general exponents (or entries in $B$). For instance, let us consider the following possibility:
\begin{equation}
    \label{qpex11}
    B = \left( \begin{array}{cc} 1 & \epsilon \\ \delta & 1 \end{array} \right)
    \:\: , \:\:\: \epsilon < 1/2 \: , \:\: \delta < 1/2 \: , \:\: \epsilon \delta \neq 1
\end{equation}
We can then find the form of the matrices $A$ and $\lambda$ such that the resulting QP system belongs to the same class of equivalence:
\begin{equation}
    \label{qpex12}
    A=B^{-1}\cdot \Gamma =
    \frac{1}{1- \epsilon \delta} \left( \begin{array}{cc}
    -1+ \epsilon /2 & -1/2 + \epsilon \\ -1/2 + \delta & -1+ \delta /2 \end{array} \right)
    \:\: , \:\:\:
    \lambda = B^{-1}\cdot \Lambda =
    \frac{ \rho }{1- \epsilon \delta} \left( \begin{array}{c}
    1- \epsilon \\ 1- \delta \end{array} \right)
\end{equation}
Notice that by construction, the class invariants of this system are $\Gamma = A_{LV2}$ and $\Lambda = \lambda _{LV2}$ displayed in (\ref{lv2inv}). Therefore Theorem 2 is applicable and the QP system (\ref{qpex11}-\ref{qpex12}) is permanent for all parameter values. Additionally, if $\rho \leq 1$ then Theorem 3 is also applicable and the unique fixed point present in the positive orthant is globally attractive in int$\{ I \!\! R^2_+ \}$.  In addition, now this system is more general than the LV of matrices $A_{LV1}$ and $\lambda _{LV1}$ in several senses. First, the QP system is reduced to the LV one in the case $\epsilon = \delta =0$ and $\rho =1/2$. Second, in the QP generalization there are more general interaction terms (i.e. quasimonomials) $\{ x_1x_2^{\epsilon},x_1^{\delta}x_2 \}$ that allow the extension of the results beyond the LV framework. Third, the QP system may describe two different conceptual situations, depending on the parameter values: if $1- \epsilon \delta >0$ the pattern is competitive, as in the original LV system:
\[
    \mbox{Pattern}(A)=
    \left( \begin{array}{cc} - & - \\ - & - \end{array} \right)
    \:\: , \:\:\:
    \mbox{Pattern}(\lambda) =
    \left( \begin{array}{c} + \\ + \end{array} \right)
\]
On the contrary, if $1- \epsilon \delta <0$ the previous two patterns are reversed and the system is cooperative, thus providing a generalization of the kind of systems for which permanence can be asserted.

\mbox{}

\mbox{}

\begin{flushleft}
{\bf 4. Global attractivity and dissipativity in dimension n}
\end{flushleft}

Some of the ideas of the previous section on global attractivity can now be extended to the framework of $n$-dimensional systems. In addition, we shall deal also with the concept of dissipativity, a treatment that complements the results on conservativity presented in
\cite{bl3}. For what is to follow we consider $I \!\! R^n_+$ as a metric space provided with the usual Euclidean distance $d(x,y)$, and define a global attractor as a set $X_0 \subset I \!\! R^n_+$ such that $\lim _{t \rightarrow \infty} d(x(t),X_0)=0$ for every initial condition $x(0) \in$ int$\{ I \!\! R^n_+ \}$. Additionally, if such global attractor is bounded the system is termed dissipative. For what is to come in this section it is also necessary to recall that a matrix $P$ is called hierarchically ordered if there exists a rearrangement of the indices such that $P_{ij} \geq 0$ whenever $i \leq j$, and additionally $P_{ii} > 0$ for every $i$. We can then state:

\mbox{}

\noindent {\bf Theorem 5.} {\em Given a $n$-dimensional QP system (\ref{qpm}) with invertible matrix $B$, if $\:\: (- \Gamma ) = - B \cdot A$ is hierarchically ordered then the system has a global attractor in $I \!\! R^n_+ $ for every initial condition $x(0) \in \mbox{int} \{ I \!\! R^n_+ \}$.}

\mbox{}

\noindent {\bf Proof.} Making use of the invertibility of $B$ a QMT of matrix $C=B^{-1}$ transforms according to (\ref{trmat}) the QP system in a $n$-dimensional LV system of matrices $A_{LV}= \Gamma$ and $\lambda_{LV}= \Lambda$, namely the QP class invariants. Since $A_{LV}= \Gamma$ the hypotheses of the theorem imply that matrix $(-A_{LV})$ is hierarchically ordered. Recall then (see \cite{hhj1}) that every $n$-dimensional discrete-time LV system for which matrix $(-A)$ is hierarchically ordered, is dissipative independently of matrix $\lambda$, so this is the case for the LV equations obtained from transformation of the QP system. Then the topological conjugacy property of QMTs \cite{bl1} allows the generalization of this result to QP systems verifying the conditions of the theorem. Note however, that in this case the full property of dissipativity cannot be mapped to class-equivalent QP systems by means of a QMT, since it is evident that the boundedness of the attracting set for the LV system is not necessarily preserved after a QMT (for instance, this may occur when some entry of $B^{-1}$ is negative). Thus, the conditions of the theorem imply the existence of a global attractor in $I \!\! R^n_+$ for the initial QP system and the result is demonstrated. $\:\:\: \Box$

\mbox{}

It is worth noting that the last result does not depend on the invariant $\Lambda$. In addition, it is clear that if it is possible to find conditions such that the boundedness of the attracting set is also preserved after the QMT in the construction of the previous proof, then dissipativity would be immediately demonstrated. A sufficient condition for this is the following one:

\mbox{}

\noindent {\bf Theorem 6.} {\em Given a $n$-dimensional QP system (\ref{qpm}) with invertible matrix $B$, if $\:\: (- \Gamma ) = - B \cdot A$ is hierarchically ordered and $B^{-1}$ is non-negative then the system is dissipative.}

\mbox{}

\noindent {\bf Proof.} Let $\{ x_i \}$ be the variables of the QP system, and $\{ y_i \}$ the variables of the canonical LV representative. According to the definition (\ref{qmt}) of QMT, and taking into account that in our case $C=B^{-1}$, both sets of variables are related by:
\[
    x_i= \prod_{j=1}^n y_j^{(B^{-1})_{ij}} \:\: , \:\:\: i = 1, \ldots , n
\]
Consequently, if $(B^{-1})_{ij} \geq 0$ for all $i,j$, the boundedness of the global attractor for the LV system is preserved after the QMT leading to the initial QP system, and then such QP system is dissipative. $\:\:\: \Box$

\mbox{}

The situation described in Theorem 6 can be applied to QP generalizations such as the one displayed in the next example:

\mbox{}

\noindent {\bf Example 2.} Consider the LV predator-prey model analyzed in \cite{lw1}:
\[
    A_{LV} = \left( \begin{array}{cc}
        -r_1 & -r_1 \mu _1 \\ r_2 \mu _2 & -r_2
    \end{array} \right) \:\:\: , \:\:\:\:
    \lambda _{LV}= \left( \begin{array}{r} r_1 \\ -r_2 \end{array} \right)
    \:\:\: , \:\:\:\:
    \{ r_1, r_2, \mu _1, \mu _2 \} >0
\]
Notice first that if Theorem 4 is applied to this system, it can be demonstrated that a unique positive fixed point exists, and it is globally asymptotically stable in $\mbox{int} \{ I \!\! R^2_+ \}$ if the following three conditions are verified: (i) $r_1 \leq 1$; (ii) $1 < \mu _2
\leq (r_1/r_2)(r_2+1)$; and (iii) $\mu_1 \mu_2 <1$. Further results not limited to specific parameter ranges can be obtained if we note that $(- \Gamma ) = -A_{LV}$ is a hierarchically ordered matrix. The system is then dissipative, according to Theorem 6. In addition, some generalizations verifying the conditions of Theorem 6 can also be constructed. For instance, more general exponents can be proposed:
\[
    B = \left( \begin{array}{cc} \rho _1 & 0 \\ 0 & \rho _2 \end{array} \right)
    \:\: , \:\:\:\: \{ \rho _1, \rho _2 \} >0
\]
For the generalized system to belong to the same class than the initial LV one, the matrices
$A$ and $\lambda$ must be:
\[
    A = B^{-1} \cdot \Gamma = \left( \begin{array}{cc}
    -r_1/ \rho _1 & -r_1 \mu _1 / \rho _1  \\ r_2 \mu _2 / \rho _2  & -r_2 / \rho _2
    \end{array} \right)
    \:\: , \:\:\:\:
    \lambda = B^{-1} \cdot \Lambda = \left( \begin{array}{r}
    r_1/ \rho _1 \\ -r_2 / \rho _2 \end{array} \right)
\]
Consequently, in spite of being a generalization (but still a predator-prey model) of the initial LV system, now the interaction exponents have been generalized while the existence of a bounded global attractor (i.e. the dissipativity property) is maintained.

%\mbox{}

%\mbox{}
\pagebreak
\begin{flushleft}
{\bf 5. Chaos criteria in dimension n}
\end{flushleft}

As it was mentioned in Section 3, the possible existence of chaotic dynamics appears naturally in connection with simpler dynamical properties such as permanence and global attractivity. It is therefore to be expected that the same analytical tools leading to the establishment of results for QP systems in previous sections provide also some new contributions in the domain of chaos. The reader is referred to \cite{diam} and \cite{maro} for the classical definitions of chaos in the sense of Diamond and in the sense of Marotto, respectively, which are adopted in what is to follow. In addition, the $n$-dimensional vector $U$ is defined as $U=(1, \ldots ,
1)^T \: \in \: I \!\! R^n$. We can now state:

\mbox{}

\noindent {\bf Theorem 7.} {\em Consider a $n$-dimensional QP system (\ref{qpm}) for which:

\noindent (a) Matrix $B$ is invertible.

\noindent (b) The class invariant $\Gamma$ is invertible.

\noindent (c) $\Gamma ^{-1} \cdot U <0$ (i.e. all the entries of $\; \Gamma ^{-1} \cdot U$ are negative).

\noindent (d) The class invariant $\Lambda$ is of the form $\Lambda = \rho \cdot U$, with $\rho \in I \! \! R$, $\rho >0$.

\noindent Then the QP system has chaotic dynamics:

\noindent (i) In the sense of Diamond if $\rho \geq 3.13$.

\noindent (ii) In the sense of Marotto if $\rho \geq 2.89$.
}

\mbox{}

\noindent {\bf Proof.} To demonstrate the result we make use of the invertibility of matrix
$B$ and the QP system is first reduced to LV form by means of a QMT (\ref{qmt}) of matrix $C=B^{-1}$. The result is a LV system with matrices $\lambda _{LV}= \Lambda = \rho \cdot U$, $\rho >0$, and $A_{LV}= \Gamma = B \cdot A$ which is invertible and verifies $A_{LV}^{-1} \cdot U<0$. For such kind of LV systems it is well-known (see \cite{doh1}) that chaotic dynamics in the sense of Diamond take place if mapping $\xi$ has a period-three orbit, where
\[
    \begin{array}{cccl}  \xi: & I \!\! R_+ & \longrightarrow & I \!\! R_+ \\
                \mbox{} & x          & \longrightarrow & \xi (x) = x \exp (\rho -x)
    \end{array}
\]
Actually such period-three orbit exists for $\rho \geq 3.13$. Similarly, it is also proved
\cite{doh1} that the LV system of matrices $A_{LV}$ and $\lambda _{LV}$ is chaotic in the sense of Marotto provided that mapping $\xi$ has a snap-back repeller, which is in fact the case for $\rho \geq 2.89$. Given that such LV system was obtained from a QP system after a QMT of matrix $C=B^{-1}$, the topological conjugacy property for QMTs \cite{bl1} allows extending these results to the initial QP system and thus implies assertions (i) and (ii) of the theorem. $\:\:\: \Box$

\mbox{}

One interesting feature of the theorem just demonstrated is that in the original LV form the result is only applicable to systems for which all the entries of matrix $\lambda$ are equal. However, now this specific form only involves the class invariant $\Lambda$, and consequently the entries of matrix $\lambda$ are different in general in the QP case, thus providing a generalization in this sense. This and other aspects of Theorem 7 are illustrated in the next example:

\mbox{}

\noindent {\bf Example 3.} We consider again the systems examined in Example 1, in which the instance provided by \cite{lw1} and some generalizations of it were investigated. Let us recall that the most general model treated there was the QP system of matrices
(\ref{qpex11}-\ref{qpex12}) with $\rho >0$. In Example 1 it was demonstrated that this system is permanent for all parameter values defined. In addition, it was also proved that if $\rho \leq 1$ then the unique positive fixed point is globally attractive in int$\{ I \!\! R^2_+ \}$. These results are interesting from the point of view of Theorem 7: according to the calculations of Example 1, notice that the class invariants of the QP system considered are:
\[
    \Gamma = \left( \begin{array}{cc} -1 & -1/2 \\ -1/2 & -1 \end{array} \right)
    \:\: , \:\:\:\:
    \Lambda = \left( \begin{array}{c} \rho \\ \rho \end{array} \right)
\]
It is then simple to check that conditions (a) to (d) of Theorem 7 are satisfied. This implies that the system is chaotic in the sense of Diamond (Marotto) if $\rho \geq 3.13$ (if $\rho \geq 2.89$). Clearly, in such chaotic scenarios the global attractivity of the fixed point is not preserved (note that the hypotheses of Theorem 3 are not verified) but permanence is still maintained (Theorem 2 remains applicable). Therefore it is demonstrated that in this family of systems permanence and chaos coexist simultaneously.

\pagebreak

\end{document}